\pdfoutput=1
\documentclass{preprint}
\usepackage[utf8]{inputenc} 
\usepackage[english]{babel}  
\usepackage[sort&compress,numbers]{natbib}

\usepackage[intlimits]{amsmath}
\usepackage{amssymb} 
\usepackage{amsthm}
\usepackage{amsfonts}
\usepackage{bm}
\usepackage{graphicx}
\usepackage{units} 

\newcommand{\h}{\mathcal{H}} 

\newcommand{\order}{\mathcal{O}} 

\newcommand{\R}{\mathbb{R}}

\newcommand\op[1]{\mathcal{#1}}
\newcommand\operate[1]{(#1)}

\DeclareMathOperator{\dx}{d}
\DeclareMathOperator{\Tr}{Tr}

\newcommand\XXint[3]{{\setbox0=\hbox{$#1{#2#3}{\int}$}
    \vcenter{\hbox{$#2#3$}}\kern-.5\wd0}}

\renewcommand{\div}{\nabla \cdot}
\newcommand{\grad}{\nabla}
\newcommand{\curl}{\nabla \times}

\newcommand{\dTr}{\op{T}_y} 

\newcommand\vek[1]{\mathbf{#1}}
\newcommand\veks[1]{\boldsymbol{#1}}
\newcommand\mat[1]{\mathbf{#1}}
\newcommand\tens[2]{\mathsf{#1}}

\newcommand{\trans}{\intercal} 

\newcommand\ut{u} 
\newcommand\trt{t} 

\newcommand\ct{T} 
\newcommand\bct{T_{I\times J}} 

\newcommand\fund[3]{\tens{#1}{#2}{\ofpt{#3}}}

\newcommand\dgamma[1]{\dx \mathrm{s}_{\pt{#1}}}

\newcommand\be{\tau}

\newcommand\comma{,\,}
\newcommand\und{\,\text{and}\,}

\newcommand\fromto[2]{\{ #1, \dots, #2 \}}

\newcommand\uu{r} 

\newcommand\pt[1]{\veks{#1}}
\newcommand\ofpt[1]{(\pt{#1})}

\newcommand\eps{\varepsilon}

\newcommand\mychi{\mathcal{X}}

\newcommand{\tikzfig}[5]{
  \begin{figure}[ht]
    \centering
    \includegraphics{#5}
    \caption{#3}
    \label{fig:#4}
  \end{figure}
}

\newcommand{\figref}[1]{Figure~\ref{#1}}
\newcommand{\Figref}[1]{Figure~\ref{#1}}

\title{\textsc{Isogeometric Boundary Element Method with Hierarchical Matrices}}

\author{ 
  J. Zechner$^{*}$, B. Marussig$^{*}$, G. Beer$^{*,\dag}$, C. Dünser$^{*}$ and %
  T. P. Fries$^{*}$%
} %

\heading{J. Zechner, B. Marussig, G. Beer, C. Dünser and T. P. Fries}

\address{$^{*}$ Institute for Structural Analysis\\%
  Graz University of Technology\\%
  Lessingstra{\ss}e 25, 8010 Graz, Austria\\%
  e-mail: ifb@tugraz.at, web-page: http://www.ifb.tugraz.at \and
  $^{\dag}$ Centre for Geotechnical and Materials Modelling\\%
  University of Newcastle\\%
  Callaghan, NSW 2308, Australia\\%
  e-mail: gernot.beer@tugraz.at
}

\keywords{isogeometric analysis, boundary element method,
  hierarchical matrices, elasticity, NURBS}

\abstract{%
  In this work we address the complexity problem of the isogeometric
  Boundary Element Method by proposing a collocation scheme for
  practical problems in linear elasticity and the application of
  hierarchical matrices. For mixed boundary value problems, a block
  system of matrices -- similar to Galerkin formulations -- is constructed
  allowing an effective application of that matrix format.
  We introduce a strategy for the geometric bisection of surfaces
  based on NURBS patches. The approximation of system matrices is
  carried out by means of kernel interpolation. Numerical results are
  shown that prove the success of the formulation.
}

\begin{document}

\section{Introduction}

In the emerging field of isogeometric analysis, Boundary Element
Methods (BEM) have gained increasing interest. This is, because for
analysis only surface descriptions are required - and Computer Aided
Geometric Design (CAD) models are based on such a boundary
description. Hence, with this combination the task of domain
discretization may be completely avoided. Still, this comes at a
prize: the numerical effort of setting up and solving the system of
equations is computationally intensive, because the system matrices
are fully populated.

Over the last decades much effort has been spent to overcome this
barrier. In context of boundary integral techniques, the fast
multipole method (FMM) \cite{rokhlin1985}, hierarchical matrices
($\h$-matrices) \cite{hackbusch1999}, the wavelet method
\cite{beylkin1991} and fast Fourier transformation based methods
\cite{phillips1994} reduce the asymptotic numerical complexity
significantly, to (almost) linear behavior.

With respect to the analysis with BEM on CAD-surfaces, early reports
on the usage of non-uniform rational B-splines (NURBS) have been
reported in \cite{Valle1994a,Rivas1996a} in the context of electric
field equations. In the field of isogeometric analysis, the strategy
was applied in \cite{Simpson2012a,Beer2013a} to practical problems of
elasticity in two dimensions and in \cite{Scott2013a,Li2011a} to
three dimensions. However, there are only few reports
\cite{Harbrecht2010a,Harbrecht2013a,Takahashi2012a} of a successful
application of fast boundary element techniques in the context of
isogeometric analysis.

In this work we present the application of the concept of
$\h$-matrices to an isogeometric NURBS-based BEM formulation for
problems in elasticity. For the geometric bisection we utilize
NURBS-features like knot insertion and the convex hull property. The
approximation of far-field matrix blocks is carried out by means of
kernel interpolation \cite{hackbusch2002}.

\section{\textsc{Isogeometric Boundary Element Method}}
\label{sec:BEM}
We consider a fixed elastic body subject to external loading.
The elastic behavior in terms of displacements $\ut$ is described by
partial differential equation
\begin{equation}
  \op{L}\ut\ofpt{x} = -\left( \lambda + 2 \mu \right) 
  \div \grad \ut\ofpt{x} + \mu \curl (\curl \ut\ofpt{x}) =
  0
\end{equation}
where $\op{L}$ denotes the Lam\'e-Navier operator. For convenience,
the boundary trace 
\begin{align}
  \label{eq:trace}
  \Tr \ut\ofpt{x} &= \lim_{\vek{x}\rightarrow\vek{y}} \ut\ofpt{x} =
  \ut\ofpt{y} & & \pt{x} \in \Omega,\, \pt{y} \in \Gamma
\end{align}
and the conormal derivative
\begin{align}
  \label{eq:conormal-derivative}
  \dTr \ut\ofpt{x} & = \lambda \div \ut\ofpt{y} \vek{n}\ofpt{y} +
  2\mu \grad \ut\ofpt{y} \cdot \vek{n}\ofpt{y} + \mu \vek{n}\ofpt{y}
  \times (\curl \ut\ofpt{y}) & & \pt{x} \in \Omega,\, \pt{y} \in
  \Gamma
\end{align}
are introduced. 
The normal $\vek{n}$ is defined to always point out of the considered
domain. The operator $\Tr$ maps displacements $\ut\ofpt{x}$ to
boundary displacements $\ut\ofpt{y}$. Involving the material law, the
conormal derivative $\dTr$ maps $\ut\ofpt{x}$ to surface traction
$\trt\ofpt{y}$. The boundary can be split into a Neumann and a
Dirichlet part such that ${\Gamma = \Gamma_N \cup \Gamma_D}$ and
${\Gamma_N \cap \Gamma_D = 0}$. This leads to the following boundary
value problem (BVP): Find a displacement field $\ut\ofpt{x}$ so that
\begin{equation}
\label{eq:BVP}
  \begin{aligned}
    \op{L} \ut\ofpt{x} & = 0 
    & & \forall \pt{x} \in \Omega\\
    \dTr\ut\ofpt{x} & = \trt\ofpt{y}  = {g}_N\ofpt{y}
    & & \forall \pt{y} \in \Gamma_N\\
    \Tr\ut\ofpt{x} &= \ut\ofpt{y} = g_D\ofpt{y} 
    & & \forall \pt{y} \in \Gamma_D.
  \end{aligned}
\end{equation}
Here, ${g}_N$ is the prescribed Neumann data in terms of surface
tractions and ${g}_D$ represents the prescribed Dirichlet data in
terms of displacements.

\subsection{Boundary Integral Equation}
\label{sec:BIE}%
The BVP \eqref{eq:BVP} can be stated in terms of an boundary integral
equation
\begin{align}
  \label{eq:BIE}
  \left( \op{C} + \op{K} \right) \ut\ofpt{x} & = \op{V}\trt\ofpt{x}
  && \forall \pt{x} \in \Gamma
\end{align}
with the weakly singular single layer operator
\begin{align}
  \label{eq:SL}
  \operate{\op{V} \trt}\ofpt{x} &= \int_\Gamma \fund{U}{2}{x,y} \trt\ofpt{y}
  \dgamma{y} & & \forall \pt{x},\pt{y} \in \Gamma
\end{align}
and the strongly singular double layer operator
\begin{align}
  \label{eq:DL}
  \operate{\op{K} \ut}\ofpt{x} &= \int_{\Gamma} \fund{T}{2}{x,y}
  \ut\ofpt{y} \dgamma{y} & & \forall \pt{x},\pt{y} \in \Gamma
  \setminus B_\eps\ofpt{x} .
\end{align}
In case of elasto-static problems $\fund{U}{2}{x,y}$ is
\emph{Kelvin's} fundamental solution for displacements and
$\fund{T}{2}{x,y} = \dTr\fund{U}{2}{x,y}$ that for tractions
\cite{beer2008}. In \eqref{eq:DL} the integral only exists as a
\emph{Cauchy principal value}, where the radius $r_{\eps}$ of a sphere
$B_\eps$ around $\pt{x}$ is treated in a limiting process
$r_{\eps}\rightarrow 0$. The remainder of that process is an integral
free term which is
\begin{align}
  \op{C} \ut\ofpt{x} &= c \ut\ofpt{x} && \forall \pt{x} \in \Gamma
\end{align}%
with $c=\nicefrac{1}{2}$ on smooth surfaces.

\subsection{Discretization with NURBS}
\label{sec:IGABEM}
In the context of isogeometric Boundary Element analysis, the geometry
is discretised by NURBS-patches
\begin{align}
  \Gamma = \Gamma_h = \bigcup_{e=1}^{E} \be_e
\end{align}
which are, in case of three dimensions ($d=3$), surface patches. Note
the equal sign for the geometry description $\Gamma$ and its
discretization $\Gamma_h$ as a unique feature: the geometry error is
zero and thus the subscript is dropped for the remainder of the
text. The function
\begin{equation}
  \label{eq:geometry-mapping}
  \mychi_\be(\pt{\uu}): \R^{d-1} \mapsto \R^d
\end{equation}
is a coordinate transformation mapping local coordinates
$\pt{\uu}=(r_1,\dots,r_{d-1})^\trans$ of the reference NURBS patch to
the global coordinates $\pt{x}=(x_1,\dots,x_{d})^\trans$ in the
Cartesian system.

B-splines form the basis of a mathematical description of the mapping
\eqref{eq:geometry-mapping} by means of NURBS. Univariate B-splines
are described by a knot vector~$\varXi=\fromto{\uu_0}{\uu_{i+p+1}}$,
which is a non-decreasing sequence of coordinates in the parametric
space, and recursively defined basis functions
\begin{equation}
\label{eq:basisfunctions}
  N_{i,p}(\uu)  = \frac{\uu-\uu_{i}}{\uu_{i+p}-\uu_{i}} \: N_{i,p-1}(\uu) 
  + \frac{\uu_{i+p+1}-\uu}{\uu_{i+p+1}-\uu_{i+1}} \: N_{i+1,p-1}(\uu) .
\end{equation}
Here, $p$ denotes the polynomial order of the B-spline and $i$ defines
the number of the knot span $\left[\uu_i, \uu_{i+1}\right)$. The initial
constant basis functions are
\begin{align}
  N_{i,0}(\uu) & = 
  \begin{cases}
    1 & \text{if~} \uu_{i}\leqslant \uu < \uu_{i+1}\\
    0 & \text{else.} \\
  \end{cases}
\end{align}
NURBS are piece-wise rational functions 
\begin{equation}
  \label{eq:NURBS}
  R_{i,p}(\uu)=\frac{N_{i,p}(\uu) w_{i}}{\sum_{j=0}^{n}N_{j,p}(\uu) w_{j}}
\end{equation}		
based on B-splines \eqref{eq:basisfunctions} weighted with $w_i$. The
basis functions $N_{i,p}$ and $R_{i,p}$ have local support and are
entirely defined by $p+2$ knots. Multivariate basis functions are
simply defined by tensor products of \eqref{eq:NURBS}. For surfaces
they are defined by
\begin{equation}
  \label{eq:NURBSnd}
  R_{\vek{i},\vek{j}}(\pt{\uu}) = \prod_{n=1}^{d-1} R_{\vek{i}_{n},\vek{j}_{n}}^n(\uu_n) 
\end{equation}
with multi-indices for the knot span $\vek{i}=\fromto{i_1}{i_{d-1}}$ and
for the order $\vek{j}=\fromto{p_1}{p_{d-1}}$ in each parametric direction.

Dropping the order-multi-index $\vek{j}$, the geometrical mapping
\eqref{eq:geometry-mapping} is now expressed by
\begin{equation}
  \label{eq:coordinate-transformation}
  \mychi_\be(\pt{\uu}) = \pt{x}\ofpt{\uu} = \sum_{\vek{i}} R_{\vek{i}}(\pt{\uu}) \pt{p}_{\vek{i}}
\end{equation}
in terms of NURBS functions and their corresponding control points
$\pt{p}=(p_1,\dots,p_d)^\trans$. In addition, Cauchy data is
discretised by the same methodology. Different to Lagrange type basis
functions, NURBS do not utilize the Kronecker delta property, hence
physical values $\vek{u}=(u_1,\dots,u_d)^\trans$ and
$\vek{t}=(t_1,\dots,t_d)^\trans$ are mapped to values in $\pt{p}$,
which are $\tilde{\vek{u}}$ and $\tilde{\vek{t}}$ marked by a
tilde. Hence, the discretization is given by
\begin{equation}
  \begin{aligned}
    \ut(\pt{x}(\pt{\uu})) \approx \vek{u}(\pt{\uu}) %
    &= \sum_{\vek{i}} \varphi_{\vek{i}}(\pt{\uu}) \tilde{\vek{u}}_{i} %
    & & \varphi \in S_h \\
    \trt(\pt{x}(\pt{\uu})) \approx \vek{t}(\pt{\uu}) %
    &= \sum_{\vek{i}} \psi_{\vek{i}}(\pt{\uu}) \tilde{\vek{t}}_{j} %
    & & \psi \in S_h^-
  \end{aligned}
\label{eq:Ansatz}  
\end{equation}
where $\varphi$ and $\psi$ are basis functions of type
\eqref{eq:NURBSnd} and $S_h$ denotes the space of basis functions,
which are at least $C^0$-continuous. With respect to physics, we
choose the Ansatz for the tractions to be discontinuous at edges or
corners. Hence, $S_h^-$ is the space of discontinuous basis functions
which are taken where the surface description
\eqref{eq:coordinate-transformation} exploits $C^0$-continuity.

\subsection{System of Equations}
By using collocation, the discretised boundary integral equation
\eqref{eq:BIE} is enforced at distinct points. Each of these points
are related to a basis function. The location of collocation points is
defined by the Greville abscissa \cite{Li2011a} except for basis
functions with $C^{-1}$ continuity. In that case, the collocation
points are slightly indented in order to avoid rank deficient system
matrices. By splitting the boundary into a Neumann $\Gamma_N$ and
Dirichlet part $\Gamma_D$ and by separating known from unknown Cauchy
data \eqref{eq:Ansatz}, a block system of equations
\begin{align}
  \label{eq:BIEdis}
  \begin{matrix}
    \pt{x}\in \Gamma_{D}: \\
    \pt{x}\in \Gamma_{N}: 
  \end{matrix} \quad
  \begin{pmatrix}
    \vek{V}_{DD} &  -\vek{K}_{DN} \\
    \vek{V}_{ND} & -\vek{K}_{NN}
  \end{pmatrix}
  \begin{pmatrix}
    \tilde{\vek{t}}_{D} \\ 
    \tilde{\vek{u}}_{N}
  \end{pmatrix}
  = 
  \begin{pmatrix}
    \vek{K}_{DD} & -\vek{V}_{DN} \\
    \vek{K}_{ND} & -\vek{V}_{NN}
  \end{pmatrix}
  \begin{pmatrix}
    \tilde{\vek{g}}_D \\ 
    \tilde{\vek{g}}_N
  \end{pmatrix}
\end{align}
with the discrete forms of \eqref{eq:SL} and \eqref{eq:DL} is created
(see \cite{zechner2013}). As a consequence of using NURBS, it is
possible to approximate known Cauchy data relatively coarsely and
differently to the unknown. 

The first subscript of the system matrices in \eqref{eq:BIEdis}
denotes the location of collocation point and the second the boundary
of the involved NURBS patches. The entries of the system matrices are
\begin{equation}
  \label{eq:matrix-entries}
  \begin{aligned}
    \mat{V}[i,j] = (\op{V}\psi_j)\ofpt{x_i} && \und && %
    \mat{K}[i,j] = \left((\op{C}+\op{K})\varphi_j\right)\ofpt{x_i}
  \end{aligned}
\end{equation}
for the $i$-th collocation point and the $j$-th basis
function. 
If the value of the basis function is zero at the collocation point,
the matrix entries are evaluated by means of standard Gauss
quadrature. For singular integrals regularisation schemes for
numerical integration are applied \cite{beer2008}. Once the matrix
entries are calculated and the known Cauchy values mapped to the
control points, the system of equation may be solved by a block
$LU$-factorisation or by means of a direct or iterative
Schur-complement solver \cite{zechner2012b}.

Due to the non-local fundamental solution $\fund{U}{2}{x,y}$ the
system matrices are fully populated so that the numerical effort for
storage and the matrix-vector-product is $\order(n^2)$. To overcome
this non-optimal complexity we apply the concept of $\h$-matrices to
\eqref{eq:BIEdis}. In the context of NURBS functions, this is
explained in the following section.

\section{Hierarchical Matrices}
\label{sec:HMAT}
In terms of the described isogeometric BEM formulation, different
approximation errors have been introduced. Firstly, by the
approximations introduced by discretization of \eqref{eq:BIE}, where
the residual is minimized in a finite number of collocation points,
and by the errors introduced evaluating integrals
\eqref{eq:matrix-entries} numerically. Secondly, the approximation of
the Cauchy data \eqref{eq:Ansatz}. Finally, the residual of iterative
solver is allowed to have a certain tolerance. Consequently, it is
reasonable to approximate the system of equations \eqref{eq:BIEdis}
itself with a similar magnitude of error. This motivated the
development of the $\h$-matrix technique by
\citet{hackbusch1999}. This matrix format provides linear complexity
up to a logarithmic factor $\order(n \log^\alpha n)$ in terms of
storage and matrix operations. For isogeometric problems of reasonable
sizes the logarithmic term is acceptable.

The matrix approximation is based on the fact, that for asymptotically
smooth integral kernels matrix blocks of well separated variables
$\pt{x}$ and $\pt{y}$ have low rank. Therefor, a partition of the
system matrices with respect to the geometry is needed. That is,
indices of matrix rows $i\in I$ and columns $j\in J$ are resorted such
that their offset corresponds somehow to their geometric
distance. Naturally, the splitting is done block-wise and categorised
into \emph{near field} and \emph{far field}. For the latter type the
variables are far away from each other and hence, the matrix block is
a candidate for approximation. 

\subsection{Geometric Bisection}
\label{sec:Bisection}
Almost every fast summation method deploys a tree to represent the
partition of matrices with general structure. The cluster tree in
context of $\h$-matrices is a binary tree and created by splitting the
geometry recursively.

\tikzfig{pics/boundingboxes}{0.0}{(a) Characteristic points with local
  bounding boxes for collcation points $Q_i$ and the support of a linear NURBS
  function $Q_j$ and (b) general binary-tree structure of a cluster
  tree $\ct$}{clustering}{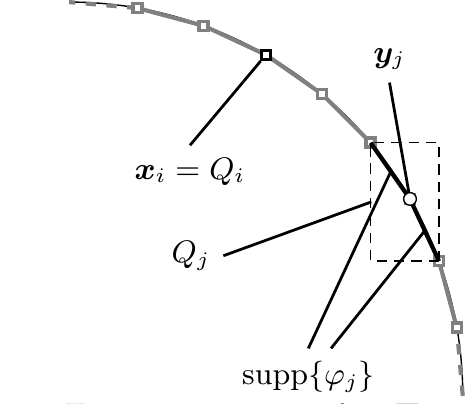}%
As shown in \figref{fig:clustering}(a), the indices $i$ and $j$ are
assigned to characteristic points $\pt{x}_i$ and $\pt{y}_j$ with local
axis parallel bounding boxes $Q_i$ and $Q_j$. Row indices of the
system matrices in \eqref{eq:BIEdis} correspond to collocation
points. Therefore $Q_i$ reduces to the characteristic point. In case
of column indices, $Q_j$ defines a bounding box around the support of
the NURBS basis function. All indices are collected to the index sets
$I$ and $J$. A cluster is the union of one or more indices of a set
including additional information stored in a label. For each set, a
labeled binary cluster tree $\ct$ is constructed. The nodes of the
tree are clusters where $t_0^0$ denotes the root cluster and is
labeled by all indices i.e. $I$, their associated positions $\pt{x}_i$
and their bounding boxes $Q_i$. Furthermore, a cluster bounding box
$B_t^\ell$ is created out of all $Q_i$ which is then geometrically
split once: $t_0^0$ gets exactly two children - the clusters $t_1^1$
and $t_2^1$. The superscript denotes the level $\ell$ in $\ct$. The
splitting is continued recursively until a stopping criterion
\begin{equation}
  \label{eq:min_leafsize}
  \mathrm{size}( t ) = \# t \leq n_{min}
\end{equation}
is fulfilled which is characterized by the minimum leaf size $n_{min}$
denoting the minimal amount of indices in a cluster. In
\figref{fig:clustering}(b) the general structure of a binary tree with
clusters is shown. If a cluster does not have any child, it is called
a leaf. In that example this is the case i.e. for $t_3^3$. The same
procedure is applied to column indices $j$ resulting to clusters $s$
and a cluster tree $\ct_J$.

For different clustering strategies we refer to the textbook of
\citet{hackbusch2009}. In the context of this work, it is suitable to
use geometrically balanced clustering. Contrary to clustering
techniques in FMM, the overall bounding box of a cluster is shrunk to
the minimum possible size with respect to the geometry $Q_i$ and $Q_j$
of the cluster-indices.
However, to perform the clustering for column indices $J$ a bounding
box $Q_j$ for each support of the NURBS functions needs to be
constructed. This is done by means of B\'{e}zier extraction and the
convex hull property. 

In our approach we generate an accumulated
knot vector $\varXi_H=\varXi_u\cup\varXi_t$ which is determined by the
individual approximation of the fields $\ut$ and $\trt$.
For a cubic curve, the following process is depicted in
\figref{fig:hull} exemplary. A B\'{e}zier extraction is performed by
means of knot insertions in $\varXi_H$ until $C^0$-continuity is
reached. The resulting control points (blue) represent a convex hull
of the NURBS curve. Hence, for each basis function $\varphi$ or $\psi$
a bounding box $Q$ of their support is generated easily by taking
these control points. For instance, the dashed box in
\figref{fig:hull} depicts the $Q_1$ for the first basis function of
the description of $\trt$ or $\ut$.
\tikzfig{pics/hull}{1.0}{B\'{e}zier extraction (blue) of a cubic NURBS
  curve described by the accumulated knot vector $\varXi_H$. The
  dashed box $Q_1$ denotes the bounding box of the support for the
  first cubic NURBS-function $R_{1,3}$ (red)}{hull}{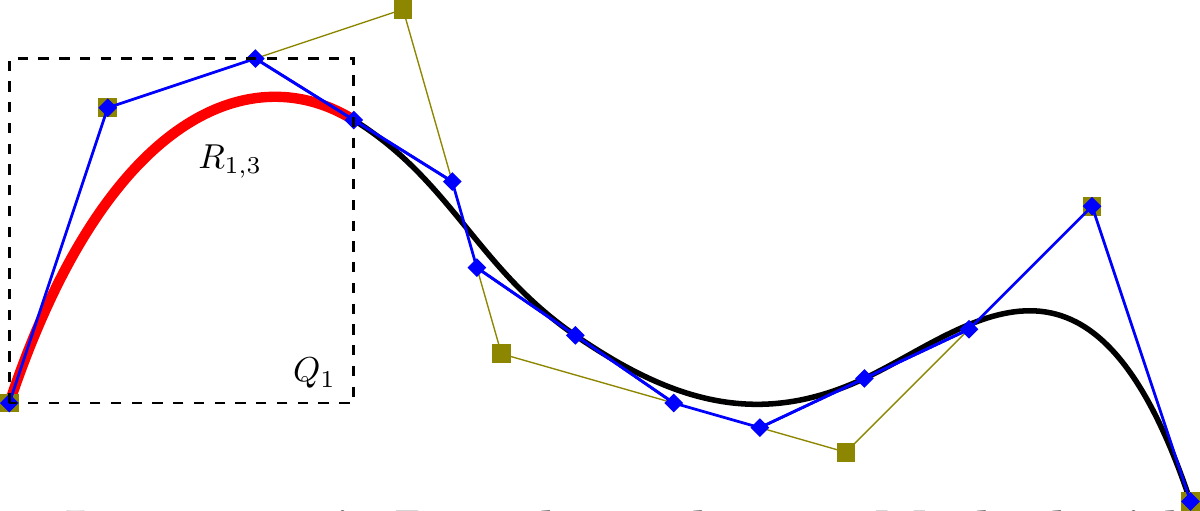}%

The structure of a $\h$-matrix is then defined by the block cluster
$\bct$ and its nodes $b = t \times s$. These nodes are constructed for
each $t$ and $s$ in the same level where an admissibility condition
\begin{equation}
  \label{eq:admissibility-condition}
  \min ( diam(B_t), diam(B_s) ) \leq \eta dist(B_t,B_s)
\end{equation}
is determined and stored. If \eqref{eq:admissibility-condition} is
fulfilled, the corresponding matrix block $\mat{M}_b$ is related to
the far field and therefor, a candidate for approximation. The block
cluster tree is now a quad tree and the basis for the partitioned
$\h$-matrix.
An example for the level-wise definition of the matrix structure is
depicted in \figref{fig:blockcluster}. Here, green matrix blocks
denote the far field. For red matrix blocks the level in $\bct$ is
increased as long as the leaf level in $t$ or $s$ is reached. Finally,
the remaining red blocks not fulfilling
\eqref{eq:admissibility-condition} define the far field. Near field
matrix blocks are evaluated with standard BEM techniques whereas far
field matrix blocks are subject to approximation. One possibility for
that is explained in the upcoming section.
%
\tikzfig{pics/blockcluster}{0.45}{Matrix partition into blocks defined
  by the block cluster tree $\bct$ in up to level $\ell=3$}{blockcluster}{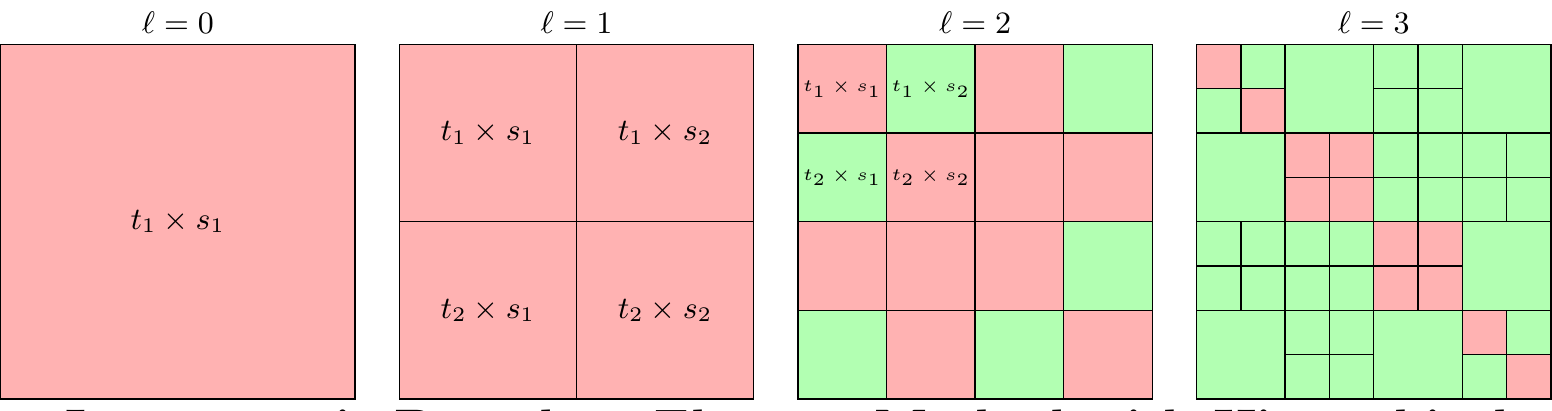}

\subsection{Matrix Approximation}
\label{sec:MatApprox}
Since the fundamental solution $\fund{U}{2}{x,y}$ is asymptotically
smooth, it is possible to separate the variables $\pt{x}$ and $\pt{y}$
to approximate the integrals \eqref{eq:SL} and \eqref{eq:DL}. Usually,
such approximations stem from Taylor or multipole expansion as well as
spherical harmonics. To avoid higher order derivatives of the kernel
function, we use the concept of kernel interpolation introduced to
$\h$-matrices by \citet{hackbusch2002}.

The fundamental solution is now interpolated by means of Lagrange
polynomials
\begin{equation}
  \label{eq:interpolation_lagrange}
  \fund{U}{1}{x,y} \approx \sum_{\nu=1}^{k} \sum_{\mu=1}^{k} 
  L_\nu\ofpt{x} \fund{U}{1}{\bar{x}_\nu,\bar{y}_\mu} L_\mu\ofpt{y}
\end{equation}
with $k$ support points defined on each of the $d$-dimensional
bounding boxes $B_t$ for $\pt{x}$ and $B_s$ for $\pt{y}$. The
interpolation functions $L_\nu$ and $L_\mu$ are represented by the
tensor product of the Lagrange polynomials in one dimension.  To get
the best approximation quality for the integral kernel, roots of
Chebyshev polynomials of the first kind are chosen for the support
points. The interpolated kernel is then taken for the representation
of single layer operator $\op{V}$ leading to
\begin{equation}
\label{eq:interpolated_SL}
  \op{V} \trt\ofpt{x} \approx \sum_{\nu=1}^{k} \sum_{\mu=1}^{k} 
  L_\nu\ofpt{x} \fund{U}{1}{\bar{x}_\nu,\bar{y}_\mu} 
  \int_\Gamma L_\mu\ofpt{y} \trt\ofpt{y} \dgamma{y}
  .
\end{equation}
As a consequence, the boundary integral in \eqref{eq:interpolated_SL}
depends only on $\pt{y}$ and is determined by Lagrange polynomials and
the traction representation. After discretization, the resulting low
rank approximation of an admissible matrix block $\mat{M}_b$ is given
by its outer product form
\begin{align}
  \label{eq:hmat:rk-representation_interpolation}
  \mat{M}_b \approx \mat{R}_k = \mat{A} \cdot \mat{S} \cdot \mat{B}^T
  && \mat{A}\in\R^{r\times k},\, \vek{S}\in\R^{k\times k},\,
  \mat{B}\in\R^{c\times k} .
\end{align}
The number of support points $k$ denote the rank of the matrices of
which the entries are given by
\begin{equation}
  \mat{A}[i,\nu] = L_\nu\ofpt{x_i} \comma\; %
  \mat{S}[\nu,\mu] = \fund{U_\ell}{1}{\bar{x}_\nu,\bar{y}_\mu}
  \;\und\; %
  \mat{B}[j,\mu] = \int_{\Gamma_e} L_\mu\ofpt{y} \varphi_j\ofpt{y} \dgamma{y}.
\end{equation}
Contrary to the quadratic storage requirement $r c$ of $\mat{M}_b$,
the requirements for $\mat{R}_k$ are only $k(r+c+k)$ which is much
smaller if $k\ll \min(r,c)$. Similar holds for the numerical effort of
a matrix-vector product. This property is the key point for the
overall reduced complexity of $\h$-matrices.

Special care has to be taken if the integral kernel depends on normal
derivatives like the fundamental solution
$\fund{T}{2}{x,y}=\dTr\fund{U}{2}{x,y}$ for the double layer operator
\eqref{eq:DL}. In that case, the conormal derivative
\eqref{eq:conormal-derivative} is shifted to the Lagrange
polynomial. The interpolated double layer potential becomes
\begin{equation*}
   \op{K} \ut\ofpt{x} \approx \sum_{\nu=1}^{k} \sum_{\mu=1}^{k} 
  L_\nu\ofpt{x} \fund{U}{1}{\bar{x}_\nu,\bar{y}_\mu} 
  \int_\Gamma \dTr L_\mu\ofpt{y} \ut\ofpt{y} \dgamma{y}
  .
\end{equation*}
It is remarkable that for both, the discrete single layer and double
layer potentials, the kernel evaluations and evaluation of $L_\nu$ at
the collocation points $\pt{x_i}$ stay the same. So do the matrices
$\mat{A}$ and $\mat{S}$. The matrix $\mat{B}$ is now defined by
\begin{align}
  \mat{B}[j,\mu] & = \int_{\Gamma_e} \dTr L_\mu\ofpt{y} \varphi_j\ofpt{y}
  \dgamma{y}.
\end{align}
For the Laplace problem $\dTr L_\mu\ofpt{y}= \grad L_\mu\ofpt{y} \cdot
\vek{n}$ holds but it can be envisaged that the implementation of the
conormal derivative for elastostatic problems
\eqref{eq:conormal-derivative} is not a straightforward task. Details on
the traction operator applied to the Lagrange polynomials are given in
the appendix of~\cite{zechner2012b}.

Since $k$ is typically chosen by the user in order to fulfil the
approximation quality, the rank of $\mat{R}_k$ might not be
optimal. In order to further reduce the storage requirement, the
matrix block is compressed by means of $QR$ decomposition. The
procedure is described in \cite{grasedyck2005}.

\section{Numerical Results}
\label{sec:Results}
To show the practicability of the described isogeometric fast boundary
element method, a numerical example in two dimensions is
presented. The approximation quality of the discretised single
\eqref{eq:SL} and double layer operator \eqref{eq:DL} is tested on a
tunnel geometry such as used in \cite{Beer2013a}. As test setting, we
chose several source points outside the domain and apply Kelvin's
fundamental solution from that points to the surface as boundary
condition. The approximation quality is measured at multiple points
inside the domain by means of the maximum norm
$\|\bullet\|_{\infty}$. \Figref{fig:convergence} shows the optimal
convergence of the described BEM formulation. As depicted in
\figref{fig:compression}, matrix compression
{$c_{\h}=\frac{\text{storage}(\mat{M})}{\text{storage}(\mat{M}_{\h})}$}
with almost linear rate is observed while accuracy is still
maintained according to the chosen interpolation quality.%
\tikzfig{pics/result_tunnel2d_error}{0.8}{Convergence of $\mat{V}$ and
  $\mat{K}$ for NURBS basis functions of order $p$ and $6$-th-order
  kernel interpolation with Lagrange polynomials}{convergence}{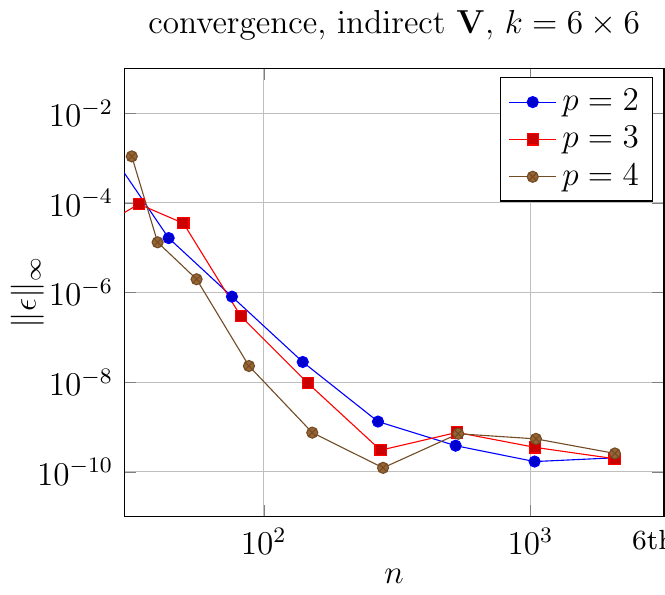}%
\tikzfig{pics/result_tunnel2d_compression}{0.8}{Compression rate of
  $\mat{V}$ and $\mat{K}$ for NURBS basis functions of order $p$ and
  $6$-th-order kernel interpolation with Lagrange
  polynomials}{compression}{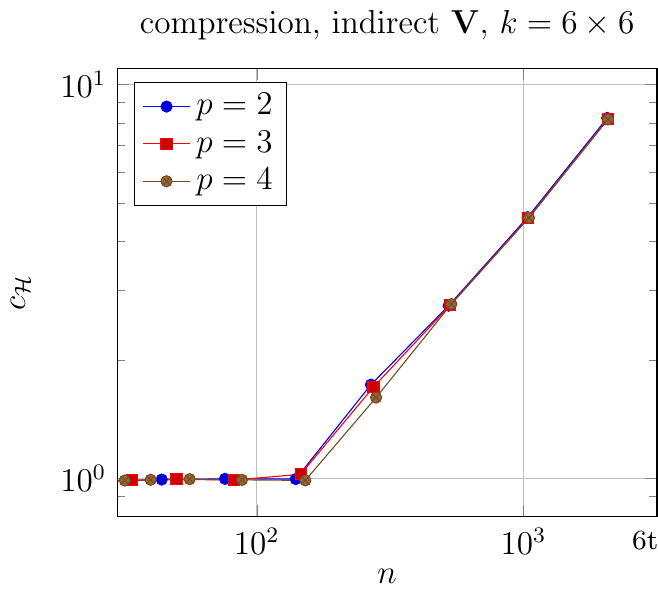}%
%

\section{\textsc{Conclusion}}
\label{sec:Conclusion}
In this work we have shown the application of the concept of
$\h$-matrices to a NURBS based, isogeometric collocation BEM. The
matrix approximation stems from the interpolation of fundamental
solutions over bounding boxes of admissible pairs of indices. For the
interpolation of the double layer operator in elasticity, the
conormal derivative to the surface is used. 
For the spatial bisection, bounding boxes enclosing the support of
NURBS functions are required. We have shown an evaluation scheme based
on knot insertion and B\'{e}zier extraction

\section*{Acknowledgment}
The authors gratefully acknowledge the financial support of the
Austrian Science Fund (FWF), Grant Number P24974-N30.

\bibliographystyle{abbrvnat}
\renewcommand*{\bibname}{REFERENCES}
\bibliography{wccm2014}

\end{document}